\documentstyle[12pt]{article}
\def \C{{\Bbb C}}

\def \bO{{\bar\Omega}}

\def \R{{\Bbb R}}
\def \O{{\Omega}}

\def \l{{\lambda}}
\def \Z{{\Bbb Z}}
\def \e{{\varepsilon}}
\def \tag{{\tilde{\frak{g}}}}
\def \ag{{\frak{g}}}

\def \o{{\omega}}
\def \K{{K\"{a}hler\ }}

\def \fd{{\bullet}}

\def \k{{\kappa}}

\def \t{{\theta}}
\def \i{{\sqrt{-1}}}

\def \d{{\partial}}
\def \ta{{\tilde A}}
\def \tom{{\tilde \Omega}}

\def \td{{\tilde D}}
\def \tc{{\tilde C}}
\def \V{{\bar V}}

\def \g{{\gamma}}
\def \d{{\partial}}
\def \z{{\bar z}}
\def \a{{\alpha}}
\def \D{{\Delta}}
\def \b{{\delta}}

\def \proof{{\noindent{\it Proof.\ \ }}}
\newtheorem{Th}{THEOREM}[section]

\newtheorem{prop}[Th]{PROPOSITION}

\newtheorem{lem}[Th]{LEMMA}

\newtheorem{corollary}[Th]{COROLLARY}

\title{Logarithmic forms with twisted coefficients}
\author{Philip A. Foth ${}^1$}

\begin{document}

\maketitle
\input amssym.def

\begin{abstract} Given a compact \K manifold, we consider the complement $U$ of
a divisor with normal crossings and a unitary local system $\tag$ on it.
We consider a differential graded Lie algebra of forms with holomorphic
logarithmic singularities and vanishing residues. We construct a spectral
sequence converging to the $L_2$ cohomology $H^{\fd}(U, \tag)_{(2)}$
corresponding to the anti-holomorphic filtration of this algebra
and compute its $E_1$ term. This differential graded Lie algebra
gives an example of a formal algebra which does not obey the $dd^c$-Lemma.
\end{abstract}

\footnotetext[1]{Sloan Doctoral Dissertation Fellow}

\section{Introduction}
\setcounter{equation}{0}

Let $X$ be a compact \K manifold and let $D\subset X$ be a divisor with normal
crossings on $X$ which has smooth irreducible components. We let $U$ be the
complement: $U=X\setminus D$, and $j:U\hookrightarrow X$ be the inclusion.
 The Hironaka resolution of singularities theorem
\cite{Hir} asserts that every smooth quasi-projective variety is birational to
a manifold $U$ of this type. Let $\tag$ be a unitary local system on $U$,
let $V$ be the corresponding holomorphic vector bundle over $U$, and let
$\V$ be its canonical holomorphic extension \cite{De1}.

   In this paper we construct a differential graded Lie algebra
(DGLA) of smooth differential forms
on $U$ with coefficients in $\V$ which have logarithmic singularities along $D$
and vanishing residues. This algebra admits the following simple local
description for a local system of rank $1$. Let $x\in X$ and let $\D$ be a
polydisc around $x$ such that in local holomorphic coordinate system $\D\cap D$
is given by $z_1\cdots z_l=0$. Then over $\D$ the sections of the sheaf $\ta^q_X(\V)$
 are generated by
$(A^{q}_X\otimes\V)_{|\D}$ and by all
expressions of the form:  $${dz_{j_1}\over
z_{j_1}}\wedge\cdots\wedge {dz_{j_k}\over z_{j_k}}\wedge \tau\otimes\nu,$$ where
$\tau\in A^{q-k}_X$, $\nu\in\Gamma(\D, \V)$, $k\le q$, $1\le j_1 <\cdots <
j_k\le l$ and such that the monodromy around $D_{j_s}:=\{ z_{j_s}=0\}$
for all $1\le s\le k$ is non-trivial.
In the holomorphic context this algebra already appeared in
a work of Timmerscheidt \cite{T'}. We say that a DGLA is formal if it is
quasi-isomorphic to its cohomology algebra.

\

{\bf THEOREM \ref{Th:main} }
{\it The $d'$-cohomology of the double
complex $\ta^{p,q}_X(\V)$ is equal to $$E_1^{p,q}=H^p(X, \bO^q\otimes
j_*\tag),$$ where $\bO^q$ is the sheaf of the anti-holomorphic $q$-forms on
$X$.}

\

We notice that in general the groups $H^p(X, \bO^q\otimes j_*\tag)$ are infinite
dimensional.

One of the important ingredients in the course of our proof is an
analogue of the $\bar\d$-Poincar\'e lemma. The DGLA $\ta^{\fd}_X(\V)$ is
naturally bigraded by holomorphic and anti-holomorphic degrees. Our Lemma
\ref{lem:d'l} tells us that locally the complex $$\cdots
\stackrel{\nabla'}{\to}\ta^{p-1,q}_X(\V)\stackrel{\nabla'}{\to}
\ta^{p,q}_X(\V)\stackrel{\nabla'}{\to}\ta^{p+1,q}_X(\V)\stackrel{\nabla'}{\to}
\cdots$$ is exact.

One of the statements that we obtain (essentially using results in
\cite{T'}) is that the corresponding spectral sequence with
$$E_1^{p,q}=H^p(X, \bO^q\otimes j_*\tag)$$ abuts to $H^{p,q}(U, \tag)_{(2)}$.
Another important conclusion is that the $d'd''$-Lemma does not hold for
our DGLA $\ta^{\fd}_X(\V)$. Nevertheless, the DGLA $\ta^{\fd}_X(\V)$ is formal,
because it is a sub-DGLA of another formal DGLA $B^{\fd}(\tag)$ and has the same
cohomology. We refer to the author's paper \cite{F} for the construction
of $B^{\fd}(\tag)$ and a proof of its formality.

\

{\bf Acknowledgments.} I am tremendously indebted to Jean-Luc Brylinski for his
advice and care.

\section{Forms with logarithmic poles}
\setcounter{equation}{0}

Let $X$ be a compact \K manifold of complex dimension $d$ endowed with a \K
form $\l$ and let $D$ be a divisor on $X$ with normal crossings which can be
represented as a union of $r$ smooth irreducible complex-analytic subvarieties
of codimension $1$:  $$D=\bigcup_{i=1}^rD_i.$$ We let $G=U(N)$ and $\ag={\frak
u}(N)$. We let $U:=X\setminus D$, we let the map $j: U\hookrightarrow X$ be the
inclusion.

A {\it differential graded Lie algebra} $(L, d)$ consists of a graded Lie
algebra $$L=\oplus_{i\ge 0}L^i, \ \ [ , ]: L^i\times L^j\to L^{i+j}, $$
satisfying for $\a\in L^i, \beta\in L^j$, and  $\g\in L^k$:
$$[\a, \beta]+(-1)^{ij}[\beta, \a]=0,$$
$$(-1)^{ki}[\a,[\beta,\g]]+(-1)^{ij}[\beta,[\g,\a]]+(-1)^{jk}[\g,
[\a,\beta]]=0$$ together with a derivation $d$ of degree $1$ (also called
a differential):  $$d: L^i\to L^{i+1}, \ d\cdot d=0, \ d[\a, \beta]=[d\a,
\beta]+(-1)^i[\a, d\beta].$$

The main examples of DGLA we will deal with include $(\ta^{i}_X(\V), d'+ d'')$,
where $$\ta^i_X(\V):=\bigoplus_{p+q=i}\ta^{p,q}_X(\V).$$

The cohomology of any
DGLA is a DGLA too, considered with zero differential. We say
that two DGLAs $(L_1, d_1)$ and $(L_2, d_2)$ are {\it quasi-isomorphic} if there
exists a third DGLA $(L_3, d_3)$ and DGLA homomorphisms $i$ and $p$: $$(L_1,
d_1)\stackrel{i}{\leftarrow} (L_3, d_3) \stackrel{p}{\to} (L_2, d_2)$$ such that
both $i$ and $p$ induce isomorphisms in cohomology. It turns out that on
the set of DGLAs quasi-isomorphism is an equivalence relation.
We also say that a DGLA is
{\it formal} if it is quasi-isomorphic to its cohomology \cite{DGMS}.

In order to construct our DGLA we introduce more notation: let
$D^{(m)}$ be the subvariety consisting of points which belong to at least $m$
irreducible components of the divisor $D$, e. g. $D^{(0)}=X$ and $D^{(1)}=D$.
Let $\td^{(m)}\to D^{(m)}$ be normalization maps and let $$v^m:\td^{(m)}\to
D^{(m)}\hookrightarrow X$$ be th composition of normalizations with inclusions.
We also let $\tc^{(m)}:=(v^m)^{-1}(D^{(m+1)})$ which is a divisor with normal
crossings on $\td^{(m)}$ (possibly empty).

Given a unitary local system $\tag$ on $U$, let $(\V , \nabla)$ be the
Deligne extension of $V$. Here $\bar V$ is a holomorphic vector bundle over $X$
and $\nabla$ is a holomorphic connection with logarithmic singularities along
$D$. We refer the reader to \cite{De1} for the details. The connections are
extended as usual to holomorphic forms with logarithmic poles.

There exist \cite{DelHod} higher residue maps $Res_m(\V )$ which are
morphisms of complexes: $$Res_m(\V ): \ \ \O^{\fd}_X\langle D\rangle
\otimes\V \to v^m_*(\O^{\fd}_{\td^{(m)}}\langle \tc^{(m)}\rangle\otimes {\bar
V}_m)[-m].$$ Here $\V _m$ is the unique vector subbundle of $(v^m)^*\V $
equipped with a unique holomorphic integrable connection $\nabla_m$ with
logarithmic poles along $\tc^{(m)}$ such that
$$Ker{\nabla_m}|_{[\td^{(m)}\setminus
\tc^{(m)}]}=(v_m)^{-1}((j_*\tag)|_{[\td^{(m)}\setminus\td^{(m+1)}]}).$$ This means
that $$(\V _m, \nabla_m)$$ is the canonical extension of
$$(v_m)^{-1}((j_*\tag)|_{[\td^{(m)}\setminus\td^{(m+1)}]}).$$ As usual, we denote
by $\O^{\fd}_X\langle D\rangle$ the complex of sheaves of holomorphic forms on
$U$ with logarithmic poles along $D$. The fact that $Res_m(\V )$ is a
morphism of complexes means that $$Res_m(\V )\circ\nabla=\nabla_m\circ Res_m(\V
).$$ Let us define $$\tom^{\fd}_X(\V ):=Ker Res_1({\bar V});$$ it is a complex
of sheaves comprised of holomorphic forms on $X\setminus D$ with 
coefficients in $\V $ which have
logarithmic poles along $D$ and have no residues. The residue maps take values in
the part of the local system invariant under all the monodromy
transformations. So when all the eigenvalues
of the
monodromy transformations are different from $1$, the complexes
$\tom^{\fd}_X(\V)$ and $\O^{\fd}_X\langle D\rangle\otimes \V$ coincide.

The complex $\tom^{\fd}_X(\V )$ admits a different description \cite{T'}.
First, on $U$ we let $\tom^p_X(\V)_{|U}:=\O^p_U\otimes\V_{|U}$. For $x\in D$ we
pick a small polycylinder $\D$ containing $x$ such that $D\cap\D$ is given by
$z_1\cdots z_l=0$ in local holomorphic coordinates $z_1, ..., z_d$. Let $T_i$,
$1\le i\le l$, be the monodromy transformation of $\tag_{|\D\setminus D}$
around $D_i$. Let $\tom^1_X(\V)_{|\D}\subset\O^1_X\langle
D\rangle\otimes\V_{|\D}$ be generated over ${\cal O}_{\D}$ by
$\O^1_{\D}\otimes_{{\cal O}_{\D}}\V_{|\D}$ and $$\O^1_{\D}\langle D_i\rangle
\otimes_{\C}(Ker(T_i-Id))^{\perp},$$ where $(Ker(T_i-Id))^{\perp}$ is the
orthogonal complement of the local subsystem $Ker (T_i-Id)$ of
${\tag}_{\D\setminus D}$, and $1\le i\le l$. Similarly we define the groups
$\tom^p_X(\V)_{|\D}$. Those sheaves glue together
nicely to a subcomplex $(\tom^{\fd}_X(\V), \nabla)$
of the complex of sheaves of
meromorphic
differential forms on $X$ with coefficients in $\V$ with logarithmic
poles along $D$.

Let us introduce a double complex $$\ta^{p,q}_X(\V):=\Gamma(X,
\tom^p_X(\V)\otimes_{{\cal O}_X}\underline{A}^{0, q}_X) , $$ where
$\underline{A}^{0,q}_X$ is the sheaf of $C^{\infty}$ differential forms on $X$
of type $(0,q)$. We need to establish certain results analogous, in a sense, to
the $d'$-Poincar\'e Lemma.

\begin{lem} Let $\k\in\R$, $0<\k <1$.
Given a complex $C^{\infty}$ function $f$ in an
open neighbourhood $W_1$ of the closure ${\bar\D}$ of a disk $\D\subset\C$
centered at zero, there exists a $C^{\infty}$ function $g$ in an open
neighbourhood $W_2\subset W_1$ of $\bar\D$ such that
$$z{\d g\over\d z}+\k g=f$$ in $W_2$. Moreover, if $f$ is
$C^{\infty}$ or holomorphic in some additional parameters, then $g$ can be
chosen to have the same properties. \label{lem:l00} \end{lem}

\proof Denote $$P=z{\d\over \d z}+\k. $$ One of the indications that the lemma
is true is that on monomials of the form $z^a\z^b$ the operator $P$ is
bijective, since $$P{z^a\z^b\over a+\k}=z^a\z^b.$$ Also, since the closure
$\bar\D$ is compact, we may and will assume that $f$ vanishes outside a compact
subset of $\C$.

Next, we decompose $$f(z)=\sum_{n\in\Z} f_n(z),$$ where $$f_n(z)={1\over
2\pi}\int_0^{2\pi}f(e^{\i\phi}z)e^{-\i n\phi} d\phi$$ and, clearly, we have
$$f_n(e^{\i\phi}z)=e^{\i n\phi}f_n(z).$$ This series is well known to be locally
uniformly convergent together with all partial derivatives. Now we claim that
when $n\ge 0$: $$f_n(z)=e^{\i n\phi}r^n h_n(r^2)$$ and when $n\le 0$:  $$
f_n(z)=e^{\i n\phi} r^{-n}h_n(r^2), $$ where $r, \phi$ are the standard polar
coordinates, $r\ge 0$, $z=re^{\i\phi}$, and $h_n(r^2)$ is a smooth $C^{\infty}$
function.  First, it is clear that $f_n(z)=e^{\i n\phi}s_n(r)$, where $s_n(r)$
is a smooth function.

Let us represent for small $r$
$$f_n(z)=\sum_{i=\max(0,n)}^k a_iz^i\z^{i-n} + o(r^k).$$

Further we have
$$f_n(z)=e^{\i n\phi}[\sum_{i=\max(0,n)}^ka_i r^{2i-n}+o(r^k)].$$

So we have $$h_n(r)=\sum_{i=\max(0,n)}^k a_i r^{2i-n-|n|} + o(r^{k-|n|}).$$
The first term has $r$ to the even power and the second term has derivatives
with respect to $r^2$ up to the order $(k-|n|)/2$. Thus $h_n$ is a smooth
$C^{\infty}$ function in $r^2$, because all its (one-sided)
derivatives $d h_n/d (r^2)$ of all orders do exist at $r=0$.

Let us now try to solve the equation $Pg_n=f_n$. First, we do it for $n\ge 0$.
We use the coordinates $(a=r^2, \phi)$, define $\d/\d a$ and $\d/\d\phi$
accordingly and notice that $$z{\d\over \d z}=a{\d\over \d a}+{1\over
2\i}{\d\over \d\phi}. $$ We will look for a solution $g_n(z)$ also in the form
$g_n(z)=e^{\i n\phi}r^nw_n(r^2)$. Applying $p$ to $g_n(z)$ one gets
$$Pg_n(z)=((n+\k)w_n(a)+a{\d w_n\over \d a})e^{\i n\phi}r^n. $$ Comparing
$Pg_n$ and $f_n$ we see that we need to solve the equation $$a{\d w_n(a)\over \d
a}+ (n+\k)w_n(a)=h_n(a)$$ on an interval $[0, \e)$. The solution is given by $$
w_n(a)={1\over a^{n+\k}} \int_0^a b^{n+\k-1} h_n(b) db.$$ It is clear right from
the formula that the solution $w_n(a)$ is smooth.

Now, we solve $Pg_n=f_n$ for $n\le 0$. Here the situation is a bit simpler and
we have  the equation $$a{\d w_n(a)\over \d a}+ \k w_n(a)=h_n(a)$$ on an
interval $[0, \e)$, which is the same as the case $n=0$ considered above.

The last thing for us to do in this Lemma is to show that the series $\sum
g_n(z)$ converges uniformly to a smooth $C^{\infty}$ function on compact
subsets. Explicitly, for $n\ge 0$ (as we have noticed the case $n\le 0$ is the
same as $n=0$) we have:  $$|{1\over a^{n+\k}}\int_0^a b^{n+\k -1} h_n(b)db|\le
{\sup_{b\in [0,a]}|h_n(b)|\over |n+\k -1|}.$$
 Now we have $$\sum_{n\ge 0}|g_n(z)|\le \sum_{n\ge 0} r^n
|w_n(r^2)| \le \sum_{n\ge 0} r^n {\sup_{b\in [0,a]}|h_n(b)|\over |n+\k -1|} <
+\infty$$ for small $r$. 
 For the first derivative we
use integration by parts:  $${d\over da}({1\over a^{n+\k}}\int_0^a
b^{n+\k-1}h_n(b)db)=$$ $$={1\over a^{n+\k}}a^{n+\k-1}h_n(a)-(n+\k){1\over
a^{n+\k+1}}\int_0^a b^{n+\k-1}h_n(b)db=$$ $$={1\over a}h_n(a)-{n+\k\over
a^{n+\k+1}}\{ [{-1\over n+\k}b^{n+\k}h_n(b)]_0^a + {1\over n+\k}\int_0^a
b^{n+\k}{d h_n\over d b}db \}=$$ $$= -{1\over a^{n+\k +1}}\int_0^a b^{n+\k} {d
h_n\over db} db.$$  This is again bounded in norm by $${\sup_{b\in [0,a]}|{d
h_n(b)\over db}|\over |n+\k|}.$$ Therefore, the absolute sum of the first
derivatives converges as well.
Clearly, the same arguments work for all the
derivatives.   $\bigcirc$

\

We denote by $d'$ the component of type $(1,0)$ of the connection
$\nabla$.

\begin{lem} For $p>0$, $x\in X$ let $\o$ be an element of
$\ta^{p,q}_{\D}(\V)$ defined in an open neighbourhood of a polydisc $\bar \D$
containing $x$ such that $d'\o=0$.  Then there exists
$\a\in\ta^{p-1,q}_{\D}(\V)$ defined in an open neighbourhood of $\bar \D$ such
that $\o=d'\a$. \label{lem:d'l} \end{lem}

\proof We can easily reduce to the case when $q=0$,
so henceforth we assume that
all the forms we deal with belong to $\ta^{\fd, 0}_{\D}(\V)$.

The statement is clear when $x\in U$.
Locally over $\D\setminus D$ our local system $\tag$
decomposes into a direct sum of local systems of rank one (since
$\pi_1(\D\setminus D)\simeq\Z^d$ is commutative). Therefore we can assume that
we have a unitary local system of rank $1$. Let $x\in D$ and first,
let us work with the case $d=1$ when the monodromy around $z=0$ is equal to
$e^{-\i\t}$, where $0<\t <2\pi$. Thus if we have a multi-valued horizontal
section $s$: $\nabla s=0$ then the monodromy operator $T$ acts on it as
$Ts=e^{-\i\t}s$. Let us consider the uni-valued section $\mu=z^{\t/2\pi}s$
which spans the canonical extension. We have: $$ \nabla\mu={\theta\over
2\pi}{dz\over z}\otimes\mu$$ and the Leibnitz rule also imply that for any
smooth function $f(z)\in C^{\infty}(\D)$, we have $$\nabla(f\mu)=[d'f+{\t\over
2\pi}f{dz\over z}]\otimes\mu= (z{\d f\over\d z}+{\t\over 2\pi}f){dz\over
z}\otimes\mu.$$  Letting $\k:=\t/2\pi$ we are within the scope of Lemma
\ref{lem:l00} and thus we are done in this case.

Next, we shall apply standard arguments \cite{Cun} to treat the case when
$x\in D^{(d)}$ and $\D\cap D$ is given by $z_1z_2\cdots z_d=0$ and the
monodromy around $D_i$ (corresponding to $z_i=0$) is $e^{-\i\t_i}$, $0 <\t_i
<2\pi$. We denote by $P_i$ the operator $z_i{\d\over \d z_i} +{\t_i\over
2\pi}$ for $1\le i\le d$.

Let  $\o={dz_k\over z_k}\wedge\psi +\beta$ where $\psi$ and $\beta$ do not
involve ${dz_j\over z_j}$ for $j\ge k$ and we will do the proof by induction on
$k$. For $k=0$ the statement is obvious. Since $0=d'\o=-{dz_k\over z_k}\wedge
d'(\psi)+ d'(\beta)$, it is clear that all the separate coefficients of the
forms $\psi$ and $\beta$ must be killed by the operators $P_j$ for $j>k$.

Any coefficient $f$ in the explicit representation of the differential form
$\psi$ (with coefficients in $\V$) can be viewed as a $C^{\infty}$ function of
the single variable $z_k$ in an open neighbourhood of the disk
${\bar \D_1}\subset\C$ obtained as the projection of $\D$ on the $k$-th
coordinate. The function $f$ is also a $C^{\infty}$ function of the auxiliary
parameters $z_1, ..., z_{k-1}$ in the polydisk
${\bar \D}_{k-1}\subset\C^{k-1}$ obtained by
the projection of $\D$ onto the first $(k-1)$ coordinates. It follows from
Lemma \ref{lem:l00} that there exists a $C^{\infty}$ function $g$ in an open
neighbourhood of $\bar\D$ such that $P_kg=f$ and that $g$ retains all the
properties of $f$ with respect to other coordinates. Let $\O$ be the
differential form obtained from $\psi$ by replacing each coefficient $f$ by the
corresponding coefficient $g$. It follows that $$d'\O=\b+ {dz_k\over z_k}\wedge
\psi, $$ where $\b$ is a differential form which does not involve ${dz_j\over
z_j}$ for $j\ge k$. Next we let $u=\o-d'\O=\beta-\b$ and we notice that $u$ does
not involve ${dz_j\over z_j}$ for $j\ge k$, since $\beta$ and $\b$ are
such.  Besides, $$ d'u=d'\o -d'd'\O=0.$$ By the induction hypothesis there
exists a $(C^{\infty})$ differential form $w$ in an open neighbourhood of
$\bar\D$ such that $u=d'w$. If we let now $\a=w+\O$ then $d'\a=\o$.

In the general case we can ignore the components of the divisor $D$ which
correspond to the trivial monodromy and by separating variables and applying the
standard $d'$-Lemma the above arguments will still work when $\D\cap D$ is given
locally by $z_1\cdots z_l=0$ and $0\le l\le d$. $\bigcirc$

\

Let us notice that another way to prove the above Lemma is to consider the
Koszul complex $$\cdots\stackrel{\cal P}{\to} E^j
\stackrel{\cal P}{\to}E^{j+1}\stackrel{\cal P}{\to}\cdots,$$ where ${\cal P}$ is
the family of commutative operators $(P_1, ..., P_l)$ and
$$E^j:=\bigoplus_{|I|=j}{\frak A}\otimes {dz_I\over z_I}, $$ where ${\frak A}$
is the algebra of germs of $C^{\infty}$ functions at zero. Since Lemma
\ref{lem:l00} tells us that $P_j$ is invertible, the complex is exact.
Let us denote by $\bO^q$ the sheaf of anti-holomorphic $q$-forms.

\begin{Th} The $d'$-cohomology of the double
complex $\ta^{p,q}_X(\V)$ is equal to $$E_1^{p,q}=H^p(X, \bO^q\otimes
j_*\tag).$$ \label{Th:main} \end{Th}

\proof The local exactness of the rows of the initial sheet of the
spectral sequence $E_0^{p,q}=\ta^{p,q}_X(\V)$ with the differential
$d'$ allows us to identify $$E_1^{p,q}=H^p(X, Ker[d': \ta^{0,q}_X(\V)\to
\ta^{1,q}_X(\V)]).$$ Next, we notice that the operator
$$P=z{\d\over \d z}+\k$$ has no kernel in $C^{\infty}$ functions on $\C$
when,  $0< \k <1$ and the kernel consists of exactly the anti-holomorphic
functions when $\k=0$. Therefore, for a local system of rank one, the space
$Ker[d': \ta^{0,q}_X(\V)\to \ta^{1,q}_X(\V)]$ consists of the anti-holomorphic
forms with support not including the components of the divisor with non-trivial
monodromy. This shows in general, that $$Ker[d': \ta^{0,q}_X(\V)\to
\ta^{1,q}_X(\V)]=\bO^q\otimes j_*\tag, \ \ \bigcirc$$

In fact, when the set of eigenvalues of the monodromies around the irreducible
components of the divisor $D$ does not include $1$, one can use cohomology with
compact support to identify $$E_1^{p,q}=H^p_c(U, \bO^q\otimes \tag).$$

We define ${\cal L}^{p,q}(\tag)_{(2)}$ to be
the sheaf of measurable $(p,q)$-forms $\a$ on $U$ with values in $\tag$ such
that both $\a$ and $(d'+d'')\a$ are locally square-integrable. Here as usual we
have the Poincar\'e metric near $D$ and the \K metric on $U$ which are
compatible. It means that the \K metric near $x\in D$ has the same asymptotic
form as the Poincar\'e metric on $\D\cap D$, where $\D$ is a small polycylinder
centered at $x$ (cf. \cite{Zu}). Set $$L^{p,q}(\tag)_{(2)}:=\Gamma(X, {\cal
L}^{p,q}(\tag)_{(2)}).$$ We let $$H^{p,q}(U,
\tag)_{(2)}:=H^q(L^{p,\fd}(\tag)_{(2)}, d'').$$ Analogously one defines the
space $H^k(U, \tag)_{(2)}$. Let us recall the following result due to
Timmerscheidt \cite{T'}:

\begin{prop} \noindent (a)
$$H^k(U,\tag)_{(2)}\simeq\bigoplus_{p+q=k}H^{p,q}(U, \tag)_{(2)}.$$

(b) The complexes of sheaves $\tom^p_X(\V)$ and ${\cal L}^{p, \fd}(\tag)_{(2)}$
are quasi-isomorphic, thus $$H^q(X, \tom^p_X(\V))\simeq H^{p,q}(U,
\tag)_{(2)}.$$

(c) The spectral sequence $$E^{p,q}_1=H^q(X, \tom^p_X(\V))$$ degenerates at
$E_1$, abuts to $H^{p+q}(X, j_*\tag)$ with $$H^k(X,
j_*\tag)\simeq\bigoplus_{p+q=k} H^q(X, \tom^p_X(\V)).$$
\label{prop:Tim} \end{prop}

We have a double complex $(\ta^{p,q}_X(\V), d', d'')$ of differential forms on
$X$ with coefficients in $\V$ with possible holomorphic logarithmic
singularities along $D$ and without residues.  The $d''$ cohomology of this
double complex is $H^q(X, \tom^p_X(\V))$ and the associated spectral sequence
with $$E^{p,q}_1=H^q(X, \tom^p_X(\V))$$ degenerates at $E_1$ by the above
Proposition. Therefore, we have the following consequence of our Theorem
\ref{Th:main}:

\begin{corollary} There exists a spectral sequence with
$$E_1^{p,q}=H^p(X, \bO^q\otimes j_*\tag)$$ which abuts to $H^{p,q}(U,
\tag)_{(2)}$.  \end{corollary}

In fact, a simple computation shows that even for $X$ a Riemann surface of genus
one with one puncture P and a local system of rank $1$ on it with non-trivial
monodromy around the puncture, the group $$H^1(X, \bO^0\otimes
j_*\tag)=H^1_c(X\setminus P, \bO^0\otimes\tag)$$  by Serre duality is dual to
the space of the sections over $U$ of the sheaf $\bO^1\otimes\tag'$. Here
$\tag'$ stands for a local system dual to $\tag$. It means that in general the
spaces $E_1^{p,q}$ are infinite-dimensional. Thus the spectral sequence that we
constructed does not degenerate at $E_1$ unless the local system extends to $X$.
We apply Proposition 5.17 from \cite{DGMS} to see that in this case the DGLA
$(\ta^{\fd}_X(\V), d'+d'')$ given by
$$\ta^i_X(\V)=\bigoplus_{p+q=i}\ta^{p,q}_X(\V)$$ does not satisfy the
$d'd''$-Lemma.

In fact, in \cite{F} the author constructs another DGLA $B^{\fd}(\tag)$,
for which the $d'd''$-Lemma holds implying its formality,
which is important for the purposes of deformation theory.
Moreover, there is a natural DGLA inclusion $$\ta^{\fd}_X(\V)\hookrightarrow
B^{\fd}(\tag),$$ which induces an isomorphism on cohomology. This means that the
DGLA $\ta^{\fd}_X(\V)$ that we consider here is formal too despite the failure
of the $d'd''$-Lemma.

\vskip 0.3in
Department of Mathematics \\ Penn State University \\ University Park, PA 16802
\\ foth@math.psu.edu

\

\noindent{\it AMS subj. class.}: \ \ 17B70, 14F05.

\end{document}